\documentclass[journal]{IEEEtran}
\usepackage{cite}
\usepackage[cmex10]{amsmath}
\usepackage{mathtools}
\usepackage{amssymb}
\usepackage{graphicx}
\usepackage{multirow}
\usepackage{booktabs}
\renewcommand{\vec}[1]{\boldsymbol{#1}}
 \normalsize

\begin{document}

\title{Estimating the Probability of Wind Ramping Events: A Data-driven Approach }

\author{Cheng~Wang, Wei~Wei,~\IEEEmembership{Member,~IEEE}, Jianhui~Wang,~\IEEEmembership{Senior Member,~IEEE}, Feng~Qiu,~\IEEEmembership{Member,~IEEE}

\thanks{This work is supported in part by the Foundation for Innovative Research Groups of the National Natural Science Foundation of China (51321005).}

\thanks{C. Wang and W. Wei are with the Department of Electrical Engineering, Tsinghua University, 100084, Beijing, China. (Email: wei-wei04@mails.tsinghua.edu.cn).}
\thanks{J. Wang and F. Qiu are with the Argonne National Laboratory, Argonne, IL 60439, USA (e-mail: jianhui.wang@anl.gov; fqiu@anl.gov).}
}

\maketitle

\begin{abstract}
This letter proposes a data-driven method for estimating the probability of wind ramping events without exploiting the exact probability distribution function (PDF) of wind power. Actual wind data validates the proposed method.
\end{abstract}


\section{Introduction}
\IEEEPARstart{L}{arge} variations of wind power, called ramping events, make it challenging to balance the load and generation in real time. A survey on the definition of a ramping event can be found in \cite{Ramp-survey}. With the point forecast results, the operator can easily identify the movement of wind power output in two successive periods that exceeds a certain threshold, i.e., a ramping event, and then schedule adequate reserve capacity so as to mitigate its impact on system frequency. However, point wind power prediction still suffers from inaccuracy as the leading time goes longer \cite{Accuracy}. Recent study proposes to forecast the confidence interval of wind power \cite{Interval}. However, the operator can hardly determine the exact movement from the confidence intervals, see Fig.~1, where the wind farm is expected to produce more power in period 1. Suppose the movement $w^e_1-w^e_2$ of point forecast constitutes a ramping event, the movement $w^u_1-w^l_2$ is certainty a more severe ramping, nevertheless, the movement $w^l_1-w^u_2$ may even not be a ramping. Moreover, the probability distribution of wind ramping capacity $w_1-w_2$ is still unclear.

\begin{figure}[h]
\centering
\includegraphics[scale=0.45]{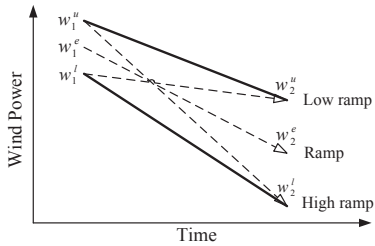}
\caption{Illustration of ramping events.}
\end{figure}

In this letter, we propose a data-driven method that can provide quantitative measure on the likelihood of ramping events given their ramping capacity, without requiring the PDF of wind power. This analysis offers statistical insights on the frequency of their occurrence and can help the operator make better generation scheduling decisions.

The exact problem studied in this letter is stated below. We have the point forecast $w^e_1$ and $w^e_2$ of wind power in two future periods, and a series of neighbouring historical data pair $(w^e_{i,1}, w^e_{i,2}), i\in \mathcal I=\{1,\dots,I\}$, in which $w^e_{i,1} \approx w^e_1$ and $w^e_{i,2} \approx w^e_2$ are met. We aim to determine the probability $\Pr[w_1-w_2 \ge R_D]$ and $\Pr[w_2-w_1 \ge R_U]$, where $R_D$ and $R_U$ are thresholds for ramp-down event and ramp-up event, which are determined by the system operator or related standards.

\section{Proposed Method}
Suppose the PDF $f(\Delta w)$ of actual wind power forecast error $\Delta w = [\Delta w_1, \Delta w_2]$ is an ambiguous multivariate function, it certainty belongs to the following functional set
\begin{equation*}
\Omega_0(\mathbb R^2) =\left\{ f(\Delta w)\left| \begin{gathered}
  f(\Delta w) \ge 0,  \forall \Delta w \in \mathbb R^2  \\
 \int_{\mathbb R^2} {f(\Delta w)\mbox{d}}w = 1
\end{gathered}  \right. \right\}   \tag{1}
\end{equation*}
The constraints in $\Omega_0$ constitute basic requirements of a PDF. Let $ \hat{f}(\Delta w)=\frac{1}{I}\sum_{i=1}^{I}(w_{i,1},w_{i,2})$ be the empirical distribution generated by the historical data, the Wasserstein ambiguity set is defined as follows
\begin{equation*}
\Omega_C =\left\{ f(\Delta w) \in \Omega_0(\mathbb R^2) \middle|
d^w_p(\mathbb P_1,\mathbb P_2) \le r \right\} \tag{2}
\end{equation*}
 where the Wasserstein metric $d^w_p(\mathbb P_1,\mathbb P_2)$ of two probability distributions $\mathbb P_1$ (described by $f(\Delta w)$) and $\mathbb P_2$ (described by $\hat f(\Delta w)$) with order $p\ge 1$ is defined by \cite{WAS_convex}
\begin{gather*}
d^w_p(\mathbb P_1,\mathbb P_2) = \inf_{\mathbb P \in \Omega_0 (\mathbb R^4)} \mathbb E_{\mathbb P}[\|z_1-z_2\|_p] \\
s.t.~  \mathbb P [z_1 \in B ] = \mathbb P_1[ z_1 \in B ], \
\forall B \in \mathcal{B} (\mathbb R^2) \\
~~~~~  \mathbb P [z_2 \in B ] = \mathbb P_2 [z_2 \in B ], \
\forall B \in \mathcal{B} (\mathbb R^2)
\end{gather*}
where $\mathcal{B} (\mathbb R^2)$ stands for all Borel sets in $\mathbb R^2$. $r$ is a measure on the distance between $\mathbb P_1$ and $\mathbb P_2$ in functional space. As $r$ tends to 0, Wasserstein ambiguity set $\Omega_C$ converges to the empirical distribution recovered from historical data.

Take the downward ramping event for example, it leads to estimate the probability $F(r_D) = \Pr[\Delta w_1 - \Delta w_2 \ge r_D]$, where $r_D = R_D - w^e_1 + w^e_2$. As the PDF $f(\Delta w)$ is not known exactly, it is prudent to investigate the worst outcome, resulting in the following optimization problem with $f(\Delta w)$ being the decision variable
\begin{equation*}
F(r_D) = \inf_{f(\Delta w) \in \Omega_{C}} \int_{\Delta w \in S_D(r_D)} f(\Delta w)\mbox{d} w      \tag{3}
\end{equation*}
where $S_D(r_D) = \{\Delta w | \Delta w_1 - \Delta w_2 \le r_D \}$. By changing the value of $r_D$, the function $F(r_D)$ provides a quasi distribution of the wind ramping capacity. It should be pointed out that for each $r_D$, the worst-case PDF $f(\Delta w)$ may not be the same.

According to Example 7 in \cite{WAS_convex}, problem (2) leads to the following convex optimization problem
\begin{gather*}
 F(r_D)  = \sup \frac{1}{I}\sum_{i=1}^{I}\beta_{i} -\gamma r  \\
s.t.\ \vec{\beta} \in \mathbb{R}^I,\ \gamma \in \mathbb{R}_{+}, \ \vec{\tau} \in \mathbb{R}_{+}^{I}  \tag{4}\\
 \beta_{i}\le 1,\  \parallel \tau_{i}\vec{s} \parallel_{q} \le \gamma, \ \forall i=1,\dots,I\\
 \beta_{i}+\tau_{i} \Delta w^h_i \le \tau_{i}r_D, \ \forall i=1,\dots,I
\end{gather*}
where $\vec{s} = [1, -1]^T$, $\Delta w_i^h = \Delta w_{i,1}-\Delta w_{i,2}$, $\Delta w_{i,1}$ and $\Delta w_{i,2}$ are the historical forecast errors,  $q$ is defined through $1/p+1/q=1$.
Some additional remarks are given.

1. Problem (3) reduces to different forms with different choice of $p$. For instance, a linear program (LP) for $p\in \{1, \infty\}$, or a second order cone program (SOCP) for $p=2$.

2. The size of $\Omega_{C}$ can be controlled through adjusting the parameter $r$.  According to \cite{Fournier_WD}, if $r$ is selected as
\begin{equation*}
r = -\log(\alpha)/I   \tag{5}
\end{equation*}
where $I$ is the number of sampled data, then the following inequality holds
\begin{equation*}
\Pr [f(\Delta w)\in \Omega_{C}] > 1-\alpha \tag{6}
\end{equation*}
where $\alpha$ is the confidence level. Equation (5) will be the main principle on the choice of $r$ in practical usage.

3. By replacing $S_D(r_D)$ in problem (3) with $S_U(r_U) = \{\Delta w | \Delta w_1 - \Delta w_2 \le r_D \}$, the convex formulation of $F(r_U)$ is similar to problem (4), except for $\Delta w_i^h=\Delta w_{i,2}-\Delta w_{i,1}$ and $\vec{s}=[-1,1]^T$.

4. The proposed method can be extended to incorporating spatial and temporal correlations, as long as there is enough historical data that produces a good reference distribution.

\section{Case Studies}

Wind data of more than 1,000 wind power plants from Jan. 1st 2004 to Jan. 2nd 2007, including the point forecasts and observed outputs, are collected from the Eastern Wind Dataset released by NREL \cite{NREL}. In this study, we select the neighbour time periods whose forecast output rests in the interval [1060 MW,1070MW], then we recover 426 data pairs. By this treatment, we have $r_D\approx R_D$ and $r_U\approx R_U$.

\begin{figure}[!t]
  \centering
  \includegraphics[scale=0.94]{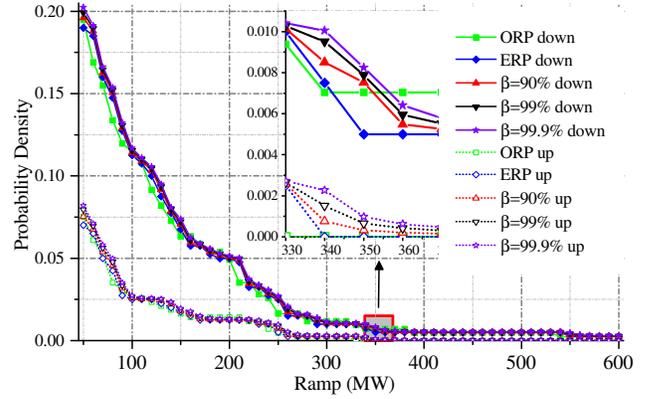}
  \caption{Quasi-PDF of ramp event when $I=400$.}
\end{figure}

In the proposed method, we choose $p=1$ such that problem (3) gives rise to an LP. We select the first 200, 300, 400 data pairs out of the 426 data pairs as samples to estimate the probability, respectively, while the empirical ramp probability (ERP) and observed ramp probability (ORP) are simulated by the known data pairs and all the data pairs, respectively. The estimated probability of downward and upward ramping event under different $R_D, R_U$ and $\alpha$ are listed in Table I, from which we can see, the probability offered by the proposed method are quite close to the real ORP and always larger than ERP. Moreover, the conservativeness can be reduced with the number of samples increasing. The quasi-PDF of ramp event when $I=400$ is shown in Fig. 2, from which we see that the ORP, ERP and estimated probability are quite close to each other. Meanwhile, the average computation time is less than 0.05 second.

\begin{table}
\tiny
\renewcommand{\arraystretch}{1.3}
  \centering
  \caption{Probability of ramping event}
  \begin{tabular}{c|c|c|c|c|c|c}
  \hline

  \hline
 \multirow{8}*{$I=200$} & $R_D (MW)$ & ORP & ERP & $\alpha=90\%$ & $\alpha=99\%$ & $\alpha=99.9\%$ \\
  \cline{2-7}
  & 200 & 0.0493 & 0.055& 0.0572 & 0.0593 & 0.0603\\
  \cline{2-7}
  & 300 & 0.0117 & 0.02& 0.0205 & 0.0210 & 0.0215\\
  \cline{2-7}
  & 400 & 0.0047 & 0.01&0.0102 & 0.0104 & 0.0106\\
  \cline{2-7}
  & $R_U (MW)$ & ORP & ERP & $\alpha=90\%$ & $\alpha=99\%$ & $\alpha=99.9\%$ \\
  \cline{2-7}
  & 200 & 0.0141 &0.015& 0.0152 & 0.0155 & 0.0157 \\
  \cline{2-7}
  & 300 & 0.0023& 0.005 & 0.0052 & 0.0054 & 0.0056 \\
  \cline{2-7}
  & 400 & 0 &0& 1.48e-4 & 2.97e-4 & 4.45e-4\\
  \hline

  \hline
  \multirow{8}*{$I=300$} & $R_D (MW)$ & ORP & ERP & $\alpha=90\%$ & $\alpha=99\%$ & $\alpha=99.9\%$ \\
  \cline{2-7}
  & 200 & 0.0493 &0.0467& 0.0481 & 0.0496 & 0.0502\\
  \cline{2-7}
  & 300 & 0.0117 & 0.0133&0.0137 & 0.014 & 0.0144\\
  \cline{2-7}
  & 400 & 0.0047 & 0.0067&0.0068 & 0.0069 & 0.007\\
  \cline{2-7}
  & $R_U (MW)$ & ORP & ERP & $\alpha=90\%$ & $\alpha=99\%$ & $\alpha=99.9\%$ \\
  \cline{2-7}
  & 200 & 0.0141 & 0.0133&0.0135 & 0.0137 & 0.0138 \\
  \cline{2-7}
  & 300 & 0.0023 & 0.0033&0.0035 & 0.0036 & 0.0037 \\
  \cline{2-7}
  & 400 & 0 &0& 9.89e-5 & 1.98e-4 & 2.97e-4\\
  \hline

  \hline
  \multirow{8}*{$I=400$} & $R_D (MW)$ & ORP & ERP& $\alpha=90\%$ & $\alpha=99\%$ & $\alpha=99.9\%$ \\
  \cline{2-7}
  & 200 & 0.0493 &0.0475& 0.0486 & 0.0497 & 0.0502\\
  \cline{2-7}
  & 300 & 0.0117 & 0.01&0.0103 & 0.0105 & 0.0108\\
  \cline{2-7}
  & 400 & 0.0047 & 0.005&0.005 & 0.0051 & 0.0051\\
  \cline{2-7}
  & $R_U (MW)$ & ORP & ERP& $\alpha=90\%$ & $\alpha=99\%$ & $\alpha=99.9\%$ \\
  \cline{2-7}
  & 200 & 0.0141 & 0.0125&0.0126 & 0.0127 & 0.0129\\
  \cline{2-7}
  & 300 & 0.0023 & 0.0025&0.0026 & 0.0027 & 0.0028\\
  \cline{2-7}
  & 400 & 0 & 0&7.42e-5 & 1.48e-4 & 2.22e-4\\
  \hline

  \hline
  \end{tabular}
\end{table}

\section{Conclusions}

A fully data-driven approach for estimating the probability of wind ramping event is proposed. Possible PDFs of the wind power forecast error are restricted in the functional Wasserstein ambiguity set. The mathematical formulation of probability estimation comes down to convex programs which are readily solvable. Case study shows that our method gives monotonically better estimation when more samples are provided.

\ifCLASSOPTIONcaptionsoff
  \newpage
\fi

\bibliographystyle{IEEEtran}
\bibliography{IEEEabrv,refs}

\end{document}